# The Key Features of the CoCalc Cloud Service Use in the Process of Mathematics and Science Teachers Training

Maiia Marienko[1[0000-0002-8087-962X]] and Kateryna Bezverbna[2[0000-0001-7088-8779]]

[1]Institute of Information Technologies and Learning Tools of NAES of Ukraine,
9, M. Berlynskoho Str., Kyiv, 04060, Ukraine
[2]Taras Shevchenko National University of Kyiv, 60, Volodymyrska St, Kyiv, 02000, Ukraine
[1]popelmaya@gmail.com
[2]katebezverbna@gmail.com

**Abstract.** The article describes the features of the evolution of the cloud-based systems in the context of mathematics and science teachers training. The notion of the systems of computer mathematics (SCM) is considered as the special kind of learning software. The lists of the special mathematical operations, functions and methods used within the structure of modern SCM are examined. Web-SCM properties are considered as general and specific. The cloud-based SCM are taken as the special stage of Web-SCM development. The criteria for an appropriate cloud service selection to support the process of mathematics disciplines learning are revealed. The key features of the CoCalc cloud service use are considered in concern to the proposed selection criteria: the available computing capacities, the diversity of tools for training and its control, the possibility of increasing the computing capability, the openness of the program code, the availability of tools for joint editing. The results of the experimental testing of CoCalc use in the process of mathematics teachers training are described. Besides, the discussion covering the issues of integration and compatibility with the other tools and services used for mathematics learning (on the example of Octave) is presented. The main trends of the CoCalc cloud service introduction into the training process and the prospects for further research are revealed.

## 1 The problem Statement

There is a need to develop new learning technologies and models of educational environment design for qualified teaching staff training, [8]. This problem can be solved using cloud technology tools. The main advantage of implementing cloud computing is the openness of access to qualitative learning resources (sometimes providing the only possible way to access the necessary resources). The idea is to identify the approaches and to assess the methods for implementing the learning components in the cloud in the process of mathematics teachers training.

**The study aims** to outline the state of the art and the current progress in the sphere of the cloud-based systems of computer mathematics development, to reveal the key features of the structure and functions of the cloud service CoCalc, relevant for its application as a tool for learning mathematics disciplines and to evaluate the promising trends of its application in the process of mathematics and science teachers training.

## 1.1 The Analysis of the Current Research

E. C. Granado, E. D. Garcia [1] have shown in their study that the Jupyter Notebook offers a rich content environment that facilitates learning while integrating the simulations performed in Python. In this way, the student learns numerically as he/she learns the concepts explained in the lesson.

The components of the author's educational mobile environment elaborated by M. A. Kislova are the following:

- the tools for organizing educational processes and mathematical calculations,
- the mobile tools to support the communication of a group of listeners,
- the tools to support the process of teaching higher mathematics and personal interaction between students and teachers [3].

Except just mobile learning the researchers have been paying much attention to the use of cloud services for learning. For example, M. A. Kislova, K. I. Slovak [2] analyzed several cloud services used in combination with the usual, traditional teaching aids in the learning process such as: Office 365, Google Apps for Education, ThinkFree Online, the advantages of these tools application for teaching mathematical disciplines were highlighted.

Substantial achievements in terms of research opportunities of the cloud technologies use in education were made by the following scientists: G. A. Aleksanyan, V. Yu. Bykov, M. Yu. Kademiya, V. M Kobysya, O. G. Kuzminskaya, V. M. Kukharenko, S. G. Lytvynova, N. V. Morse, V. S. Mkrtchyan, A. S. Sviridenko, Z. S. Seydametova, S. O. Semerikov, O. M. Spirin, L. V. Rozhdestvenskaya, Y. V. Trius, M. P. Shyshkina, B. B. Yarmakhov, M. I. Zhaldak, and others. These issues were also investigated by M. Armbrust, R. Griffith, Y. Khmelevsky, M. Miller, K. Subramanian, N. Sultan, W. Chang, P. Thomas, A. Fox, et al.

The use of the cloud services has a number of advantages (according to N. V. Rashevskaya) [6]:

- data will be available from any device with the access to the Internet;
- working with all educational materials without installing third-party software;
- the opportunity to study anywhere and anytime (outside the classroom);
- the possibility of using blended and distance learning.

There is a separate group is research studies devoted to the problems of cloud technologies use for teachers training. Among them, there are the works of T. L. Arkhipova, N. V. Bakhmat, T. V. Zaitseva, Yu. G. Lotyuk, N. V. Soroko, M. A. Shynenko et al.

The use of the Sage computer algebra system for learning mathematics was studied by P. Zimmermann, A. Casamayou, N. Cohen, G. Connan, T. Dumont, L. Fousse, F. Maltey, M. Meulien, M. Mezzarobba, C. Pernet, N. M Thiery, E. Bray, J. Cremona, M. Forets, A. Ghitza, H. Thomas [9]. The researchers were exploring the functionality of symbolic computing systems, giving examples of the use of separate commands to solve common mathematical problems.

In the study of A. B. Liecharlie [4], the CoCalc cloud service was used to process a large array of data.

M. D. Ruiz and F. Torralbo [7] investigated the advantages and disadvantages of using the CoCalc cloud service in training. The authors have described the platform, the methodological system and the peculiarities of teaching in the context of different modes of use of the cloud service. Some methods of teaching different subjects using CoCalc were exposed.

As there are new facilities and modes of using the cloud services in education due to the latest achievements and developments in this area, in particular, the current progress of CoCalc applications, the issues of innovative use of this kind of services need further attention. Among them, there are issues of the search for the new forms, methods and possible models of CoCalc use in the process of training mathematics teachers to support the teaching of mathematical disciplines.

## 2 The Main Results of the Study

### 2.1 The Main Features of Modern SCM Composition

The modern SCM can be of different types, application areas and architecture; still, they appear to have similar key features of their structure and composition [2]. They include usually the next main components:

— one of the main components is the computing core of the system;
— standard default codes for use in complex user functions and procedures;
— intuitive user interface that will allow you to perform quickly complex calculations using a standard set of functions;
— powerful graphics packages that may be used in mathematics and other areas;
— packages of additional modules that significantly expand the capabilities of SCM;
— additional modules, libraries and functions that are not included by default in the kernel;
— flexible help system that simplifies the use of certain functions, libraries, templates.

SCM provides strong support for the study of mathematical disciplines:

— the ability to open parentheses in expressions containing characters;
— the calculations of the values of numerical expressions;
— the use of symbolic values of variables to calculate a symbolic expression;
— to simplify expressions, including the opening of parentheses;
— finding the exact and approximate roots of equations or systems of equations;
— performing a number of tasks in mathematical analysis;

— construction of graphs of functions and surfaces, images of vectors on the plane and in space;
— a number of functions to perform tasks of linear algebra and others.

## 2.2 The Evolution of the Cloud-Based Systems of Computer Mathematics

In the process of SCM evolution the special kind of them such as Web-SCM emerged along with their gradual transformation into cloud-based systems. In recent years, the cloud-based versions of mathematical packages from major providers, such as Maple Net, MATLAB web-server, WebMathematica and others have been supplied. A variation of this type of systems is CoCalc, the cloud-based version of Web-SCM Sage.

SCM has undergone some changes over time. During the 60's and 90's of the XX century SCM were mostly local, ie required installation on a local computer.

In the 90s of the XX century and until the first decade of the XXI century, Web-SCM appeared, which allows you to work with SCM directly via the Internet and configure a web browser. It was the second stage of SCM evolution. The main features of such SCM are:

— the computing core does not require pre-installation on the user's device;
— all calculations take place on the web server (user standby capacities may be low);
— all calculations can be organized using a web browser.
— Web-SCM also has specific functions:
— does not depend on the technical characteristics of the user's personal device;
— you can use any browser;
— intuitive ease of introduction;
— constant access to educational materials anywhere and anytime, etc.

The most commonly used are Web-SCMs such as MapleNet, MathCAD Application Server (MAS), webMathematica, Matlab Web Server (MWS), Sage and wxMaxima. Web-SCM has an intuitive interface, powerful graphics packages for building and visualizing mathematical objects, functions and methods.

Specific characteristics of modern SCM are:

— the use of the elements of mathematical functions and programming languages;
— integration with other software;
— the capability of editing and printing of mathematical texts.

The use of Web-SCM had certain advantages: the availability of powerful tools for building and researching mathematical models for various mathematical disciplines.

The first version of SAGE emerged in 2006. The first cloud systems were developed in 2009, it caused the third stage in the development of SCM. CoCalc (new name) is a cloud service that combines all the features and tools of SAGE and Web-SCM SAGE.

### 2.3 The Criteria for Cloud Services Selection

The following selection criteria may be used for selecting a cloud service, appropriate for the certain educational needs:

- computing power provided for the work of an individual user, as well as how much data may be processed simultaneously;
- whether there are tools for organizing and monitoring the learning process (for certain disciplines, setting different levels of learning tasks, reviewing the process of implementation and evaluation of learning results, both for groups of users and individual students);
- availability of different tariff plans to expand the computing capabilities of the cloud service (availability of additional functions);
- open source code and free access to the installation or integration of additional modules and applications (including author's tools) along with the basic set of tools being protected;
- availability of tools for creating, storing, joint editing of resources or files in various formats (the most common and compatible with modern devices).

### 2.4 The Comparison of SageMathCloud and CoCalc

Sage is a software that implements mathematical algorithms in different contexts. To begin with, it can be used as a scientific pocket calculator and can manipulate all sorts of numbers, from integer and rational numbers to numerical approximations of real and complex numbers with arbitrary precision. However, mathematical calculations go far beyond numbers: Sage is a system of computer algebra; for example, it can be used to support solving linear equations or search for solutions, simplify expressions. In the process of analysis, Sage can manipulate by the expressions relating to square roots, exponents, logarithms, or trigonometric functions: integration, calculating limits, simplifying sums, series, solving certain differential equations, and more. In linear algebra operations with vectors, matrices, and subspaces can be performed. It can also help illustrate individual tasks in probability theory, statistics, and combinatory [9].

Sage service will provide the user with consistent access to functions across a wide area of mathematics - from group theory to numerical analysis, and further to two- and three-dimensional visualization, animation, networking, databases. Using all the necessary functions within the unified cloud service is convenient for a user who needn't transfer data between multiple tools and learn the syntax of several programming languages.

The specifics of using CoCalc is primarily the capability of projects creation (both group and individual) and uploading training materials to them. After that users will be able to share them. There is a chronology of user actions with individual training material and with the project as a whole. These features are useful in the process of evaluating each user individually.

CoCalc is the online computing environment [1] that provides several unique tools of learning use: it enables Jupyter Notebook file sharing, which empowers group work and supports the course and task manager, which simplifies the teacher's task of supplying support material, assigning task solving and providing comments and grades for students. It is also free of charge and therefore students have access to the materials at the end of the course, which occupies modelling and explanations that may be useful in other subjects of the following courses. It also has paid packages that offer access to more stable servers with more computer resources.

The freeware software that is installed is much diversified, covers many programming languages, and can be easily expanded by installing other programs on Linux in user space. CoCalc user personal space is organized into different projects, each with its files that can be generated on the platform or downloaded from a user's computer. On the other hand, it also has tutorials on different aspects of both platform use and programming in different languages. Regarding the course management, the teacher can add students by their email address, assign tasks, automatically form a specific folder for all course participants, keep track of the edits to these assignments, and return a corrected copy to all students with the one click of a button. Another useful feature of CoCalc as for a teacher point of view is that it allows simultaneous editing of Jupyter Notebook documents in the Google Documents style, thanks to the proprietary version of the Jupyter Notebook web interface. This add-on, not included in the default Jupyter notepad, greatly facilitates group work on the task and is complemented by integrated chat on the platform itself. All of these options make CoCalc a complete and easy-to-use solution for integrating Jupyter files into an academic course. Also, the openness of the code led to the emergence of a local version of the platform, which could be installed on any computer. This local version can be downloaded from the GitHub development page and, if the local computer is powerful enough, used for teaching a Jupyter Notebook course [1].

CoCalc (formerly called SageMathCloud) is a cloud computing and math calculations management platform. CoCalc directly supports the Sage worksheet, which interactively evaluates the Sage code. SageMath (formerly Sage or SAGE, an algebra and geometry system for numerical operations and research) is a computer algebra system that has various functions that cover many aspects of mathematics, including algebra, combinatorics, graph theory, numerical analysis, calculus, and statistics. The new name CoCalc, which stands for Collaborative Calculations, reflects the evolution of the conventional cloud computing and calculation platform. Spreadsheets in CoCalc support Markdown and HTML as design and R, Octave, Cython, Julia and others for programming other than Sage [4]. The main difference between CoCalc and SageMathCloud is that CoCalc is primarily focused on collaborating on any file in the system, editing and communicating during teamwork and organizing the students group learning process. While the main benefit of SageMathCloud was the use of the cloud computing approach. Almost every CoCalc file contains a built-in chat, participant group change history, and an easy (intuitive) feature to add new users to the project.

CoCalc is an open-source software provided by SageMath Inc. The creator and principal developer of CoCalc is William Stein, a professor of mathematics at Washington University who also created the Sage software system. The initial development of CoCalc was funded by Washington University and grants from the National Science Foundation and Google. At the moment, CoCalc is largely funded by paid users.

## 3 The Results of the Experimental Testing of the Use of CoCalc in the Educational Process

To test the effectiveness of using CoCalc as a tool for teaching mathematical disciplines for pre-service teachers of mathematics the pedagogical experiment was conducted in 2017, at the Pedagogical Institute of the Kryvyi Rig National University (Kryvyi Rig).

The experimental group and the control group were formed according to the following principle: in the experimental group there were students who studied the author's method of using CoCalc (EG), in the control group there were students who studied the traditional method (CG). The obtained results proved that the level of teachers' professional competencies formation in the experimental and control groups coincided with the level of significance $\alpha = 0.05$ before the formation stage [5].

Comparing the levels of professional competences formation in control and experimental groups it was possible to observe the increase of a share of students with average and high level of professional competence after the formation stage of the experiment.

Comparison of the distribution of the experimental and the control groups of students after the formation stage of the experiment (scale levels): high – from 16% to 20%, sufficient – from 26% to 36%, average – from 27% to 28%, low – from 31% to 16%. During the study only some components of subject, professional and practical and technological competencies have been compared at four levels (high, sufficient, medium and low), among them: subject-pedagogical, mathematical and information-technological competencies. According to Fisher's test, the data at the beginning of the experiment and at the end were compared.

The analysis of the results of the forming stage of the pedagogical experiment showed that the distribution of the levels of the formation of professional competencies in the experimental and control groups of mathematics trainee teachers has statistically significant differences due to the implementation of the developed method of using the cloud service CoCalc, which confirms the hypothesis of the study [5].

During the various teaching experiences, several shortcomings of the cloud service use were identified. The first disadvantage is the inability to prepare self-assessment questionnaires (simulators, tests). This tool, available on many training platforms such as Moodle, is not currently available on CoCalc. This tool would allow a teacher to learn the level of theoretical content knowledge of a student on a subject before the practical tasks. However, the CoCalc extension is planed [7], that will allow the teacher to evaluate the exercises by developing tests automatically. The second drawback concerns to only in the free tariff plan, which limits the calculating of large amounts of data, or the installation of additional libraries if necessary for a specific task. However, in the process of teachers training, the high computation and accuracy process maybe not so important. An alternative that addresses this second drawback is the purchase of a paid tariff plan from CoCalc. Another option is to install and maintain a Jupyter server using JupyterHub software [7], in which the SageMath character system should be installed as a Jupyter kernel. The second option will allow you to perform intensive calculations (according to server performance), as well as the possibility to install the Jupyter Nbgrader extension [7], to implement a form of self-assessment for programming tasks. One more option, if the agency has a contract with Google, is to use Google

Colaboratory, which allows to implement Jupyter and share an interface similar to Google Drive.

The third limitation that was found concerns to the assessment of group work, since there is no option to create or evaluate activities for a group of users in a whole, which also allows you to evaluate the participation of each small group within it.

## 4 Discussion

There is a way to examine the difference between working with software in local-based and cloud-based versions comparing the use of different programming tools of CoCalc cloud service (by the example of Octave).

Octave is available in CoCalc along with other kinds of languages such as Python, Java, C, and so on. Octave is an interpreted language designed for mathematical computing, it is largely compatible with MathLab. This is the reason that this tool is becoming more and more popular and useful in the process of math teachers training. In particular, in the Taras Shevchenko National University of Kyiv (Ukraine) there is a separate laboratory course devoted to it.

The use of Octave based on CoCalc enables the user to work with Octave in a cloud and also to use it in collaboration. Sagews may provide the mathematical functions (for working with matrices, vectors, graphs, mathematical constants, equations), but the syntax may be different and not compatible with MathLab. The source code for Octave is freely distributable under the terms of the GNU General Public License (GPL) it makes possible to use free all of its functions in the learning process.

In the cloud version, there are the same commands as in the local one (Fig. 1), and also there are some additional possibilities. There are tools for easy collaboration such as video chat, chat commands, code history of each user, the ability to save files (Fig. 2).

The advantage of the cloud version is that it can be accessed from different computers, even from a smartphone.

The cloud-based interface of Octave used through CoCalc may slightly differ from its local version, though not very significant, as its capabilities remain. When you create a file with the extension .m, you open an editor where you can edit the text, you can cut, copy, paste, there is a search in the text. The '> _Shell' key lets you run the code. The 'file' key enables you to work with the file (Fig.3).

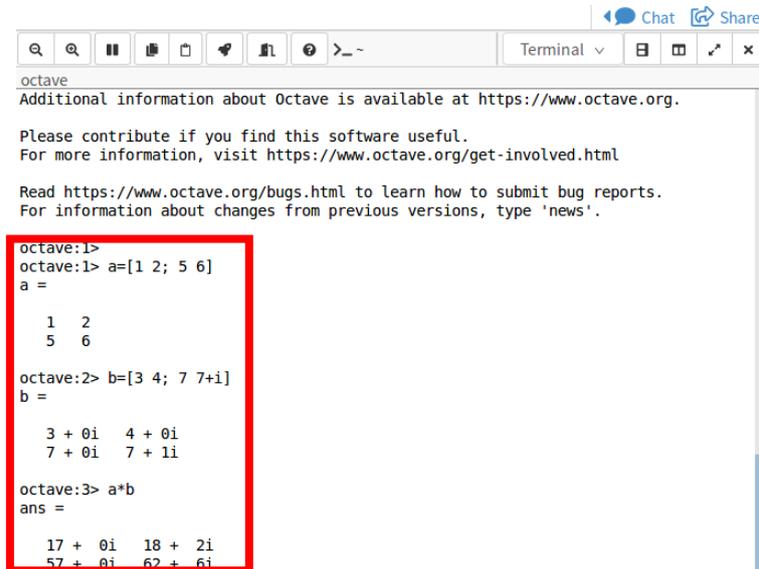

**Fig. 1.** The Octave code is running in the terminal

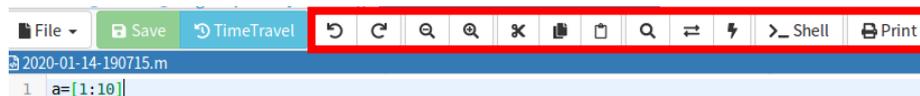

**Fig. 2.** Open editor and control keys

There are different ways to work with Octave in CoCalc. You can simply create a file with the extension .m, then you can open the terminal and run the file in it.

The .term file is used to use Linux terminal commands. In the terminal, you can run Octave or create a file with extention .m of .x11.

In the file with extention .x11 you can run different applications. There is possible to open Libre OfficeWriter, the environment for Python and so on and also use Octave as in the local version.

Now the Octave language is penetrating the math teachers' education process. This language is a useful tool for teachers training as it enables learners to perform different mathematical calculations and operate with functions. Comparison of using the local version of Octave and the cloud one (based on CoCalc) is a prospective trend for further research.

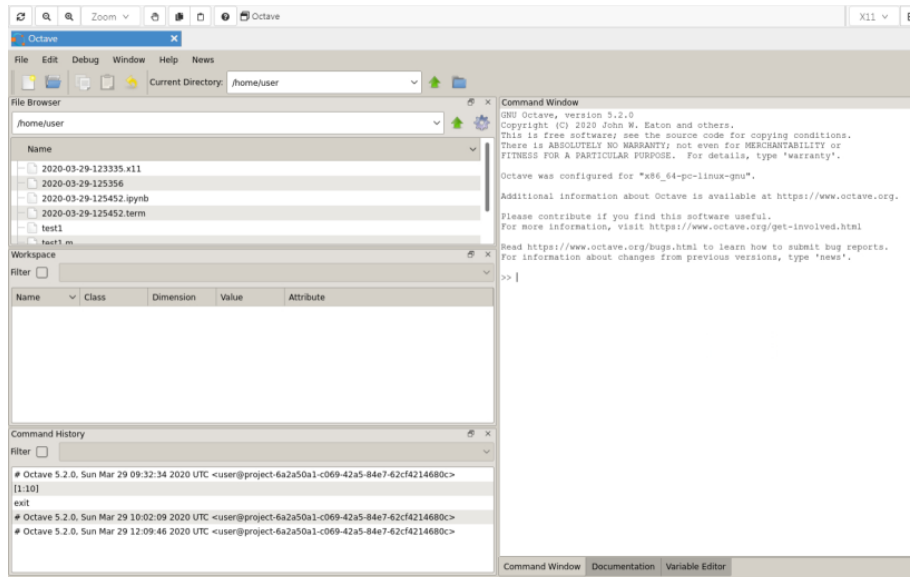

**Fig. 3.** Octave Editor: general view

## 5     Conclusions

The emergence and development of new forms of learning, focused on the collaboration of groups of students and joint activities in the online mode was caused by the availability of cloud services (including specialized). It is justified that the cloud services may be purposely used in the process of science and mathematics teachers training as tools for:

─ group communication of listeners (synchronous and asynchronous);
─ collaboration in real time;
─ processing, systematization and storage of educational data.

Among the areas of CoCalc use in the science and math teachers training process are the next:

─ organization of communication during the educational process;
─ support of various forms of educational process (individual and group) both classroom and extracurricular;
─ support for the organization of training;
─ visualization of interpretations of mathematical models, mathematical abstractions by visualization;
─ providing open access using a common interface and reliable software;
─ increasing mobility in time and space;
─ compatibility with other types of mathematical software, as well as the possibility of cloud access to various software packages, including Octave based on CoCalc;

- formation of the learning environment and filling its content during the learning process.

The described teaching experience with the use of CoCalc demonstrates its potential to become a ground for a variety of learning methods.

The teacher, through the use of the control and monitoring tools offered by the system, can supply and control students performance that is an added benefit of using this platform. Finally, this service is also suitable for use as the research collaborative platform.

## References

1. Díaz García, E., Cabrera Granado, E.: Manual de uso de Jupyter Notebook para aplicaciones docents. Universidad Complutense de Madrid. https://eprints.ucm.es/48304/ (2018). Accessed 28 March 2020
2. Kislova, M. A., Slovak, K. I.: Clouds for teaching mathematical disciplines. Modern computer technologies. 11, 53-58 (2013)
3. Kislova, M. A.: Development of a mobile educational environment for higher mathematics in the preparation of electromechanical engineers: author. diss. ... Cand. ped. Sciences: 13.00.10 – information and communication technologies in education. Ph.D. thesis, Kyiv (2015)
4. Liecharlie, A. B.: Penentuan rute optimal transportasi kontainer di kalimantan barat. Jurnal TIN Universitas Tanjungpura. **2**(1), 16-21 (2018)
5. Popel, M. V.: Using Cocalc as a Training Tool for Mathematics Teachers'pre-Service Training. Information Technologies and Learning Tools. 6(68), 251-261 (2018)
6. Rashevskaya, N.: Cloud Computing in Higher Mathematics Teaching in Technical Universities. In: Cloud Technologies in Education: Materials of the All-Ukrainian Scientific-Methodological Internet Seminar (Krivoy Rog – Kyiv – Cherkasy – Kharkiv, December 21, 2012 ), pp. 127-129. KMI Publishing Department, Kryvyi Rih (2012)
7. Ruiz, M. D., Torralbo, F.: Tres experiencias docentes con una plataforma en la nube. Enseñanza y Aprendizaje de Ingeniería de Computadores. **8**, 27-36 (2018)
8. Shyshkina M.: Holistic Approach to Training of ICT Skilled Educational Personnel. In: Ermolayev, V. (ed.) Proceedings of the 9th International Conference on ICT in Education, Research and Industrial Applications : Integration, Harmonization and Knowledge Transferine, vol. 1000, pp. 436-445. CEUR Workshop Proceedings. http://ceur-ws.org/Vol-1000/ICTERI-2013-p-436-445-MRDL.pdf (2013). Accessed 28 March 2020
9. Zimmermann, P., et al.: Computational Mathematics with SageMath. SIAM- Society for Industrial and Applied Mathematics (2018)